\documentclass[12pt, reqno]{amsart}
\usepackage{amsmath, amsthm, amscd, amsfonts, amssymb, graphicx, color}
\usepackage{setspace}
\usepackage{mathrsfs}
\usepackage{multicol}
\usepackage[bookmarksnumbered, colorlinks, plainpages]{hyperref}
\hypersetup{colorlinks=true,linkcolor=red, anchorcolor=green, citecolor=cyan, urlcolor=red, filecolor=magenta, pdftoolbar=true}

\newtheorem{theorem}{Theorem}[section]
\newtheorem{lemma}[theorem]{Lemma}
\newtheorem{proposition}[theorem]{Proposition}
\newtheorem{corollary}[theorem]{Corollary}
\theoremstyle{definition}

\theoremstyle{remark}

\numberwithin{equation}{section}

\newcommand{\CC}{\mathbb {C}}

\begin{document}
\setcounter{page}{1}
\title[Topological structures of generalized   Volterra-type integral  operators]{ Topological structures of generalized  Volterra-type integral  operators}
\author[Tesfa  Mengestie]{Tesfa  Mengestie }
\address{Department of Mathematical Sciences \\
Western Norway University of Applied Sciences\\
Klingenbergvegen 8, N-5414 Stord, Norway}
\email{Tesfa.Mengestie@hvl.no}
\author [Mafuz Worku]{ Mafuz Worku}
\address{Department of Mathematics,
Addis Ababa University, Ethiopia}
\email{mafuzhumer@gmail.com}
\begin{abstract}
We study  the generalized Volterra-type integral and composition operators acting on
 the classical Fock spaces. We first  characterize  various properties of the operators  in terms of growth and integrability  conditions  which are  simpler  to apply than  those already known Berezin type characterizations. Then, we apply these conditions to study the compact   and  Schatten $\mathcal{S}_p$ class difference  topological structures of the  space of the operators. In particular, we proved that the difference of two   Volterra-type integral operators is compact if and only if both are compact.
\end{abstract}
\maketitle
\section{introduction}
Ever since A. Aleman  and A. Siskakis published their seminal works \cite{Alsi1,Alsi2} on  integral operators acting on Hardy and Bergman spaces,  Volterra-type integral operators have been
studied on many Banach spaces over several domains. For  holomorphic functions $f$ and $g$ in a given domain, we define  the  Volterra-type
integral operator $V_g$ and its companion $J_g$   by
\begin{align*}
 V_gf(z)= \int_0^z f(w)g'(w) dw \ \ \ \text{and} \ \  \ J_gf(z)= \int_0^z f'(w)g(w) dw.
\end{align*}
Applying integration by parts in any one of the above integrals
gives the relation \begin{align}
\label{parts}V_g f+ J_g f= M_g f-f(0)g(0),\end{align} where $M_g
f= g f$ is the multiplication operator of symbol $g$.
  Studies on these  operators have been  mainly aiming to descsribe  the connection between
their operator-theoretic behaviours with the function-theoretic properties of  the inducing
symbols $g$.    For more information on the subject,
we refer to \cite{Alman,ALC,TM00} and the related references
therein.

In 2008, S. Li and S. Stevi\'c  extended $V_g$ and $J_g$ to  the
 operators
\begin{align}
\label{OKK}
V_{(g,\psi)} f(z)= \int_0^z f(\psi(w)) g'(w) dw, \ \text{and}\ \  J_{(g,\psi)} f(z)= \int_0^z f'(\psi(w)) g(w) dw,
\end{align}
and studied some of  their operator-theoretic properties in terms of  properties of the pairs   $(g,\psi)$  on some spaces of analytic
functions on the unit disk \cite{LIS,jmaa345}. For more  recent
results on this  class of  operators, one may also consult the materials for instance  in  \cite{TM, TM0, UK2}.

 We may note in passing that if we, in particular, set $\psi(z)= z,$ then  the operators $ V_{(g,\psi)}$ and $J_{(g,\psi)}$  become  respectively   the operators $V_g$ and $J_g.$  On the other hand, the choice $g= \psi'$ reduce  the operators $J_{(g,\psi)}$ to  the composition operators $C_\psi$ taking $f$ to $f(\psi)$ up to a constant.  Owing to this fact, the operators in \eqref{OKK} are often referred as the generalized  Volterra-type integral operators  and  generalized composition  operators.
 Such operators have found applications for example  in the study of linear isometries of
spaces of analytic functions. If  $S^p$ denotes the
space of all analytic functions $f$ in the unit disc for which its
derivative $f'$ belongs to the Hardy space $H^p, $  then it has been shown that   for
$p\neq2,$ any surjective isometry $T$ of $S^p$ under the norm
$\|f\|_{S^p}= |f(0)|+ \|f'\|_{H^p}$ is of the form $$Tf= \lambda
f(0)+ \lambda J_{(g,\psi)}f$$ for some unimodular $\lambda$ in $ \CC$,
 a nonconstant inner function $\psi$ and a function $g$ belonging to the space  $H^p$ \cite{FJ}.

Let $\CC$ be the complex plane and  $0<p<\infty $.  Then   the classical
Fock spaces $\mathcal{F}_p$ consist of entire functions $f$ for which
\begin{equation*}
\|f\|_{p}^p=  \frac{p}{2\pi}\int_{\CC} |f(z)|^p
e^{-\frac{p }{2}|z|^2} dA(z) <\infty,
\end{equation*} where $dA$  denotes the Lebesgue area  measure.  In particular, $\mathcal{F}_2$ is a
reproducing kernel Hilbert space with kernel and normalized reproducing kernel
functions   given by the  explicit formulas
\begin{align*} K_{w}(z)= e^{ \overline{w}z}\ \   \text{and} \ \
k_{w}(z)= e^{\overline{w}z-\frac{|w|^2}{2}}.
\end{align*}
 The kernel function $K_{w}$ belongs to all the
  Fock space $\mathcal{F}_p$  with  norms
  \begin{align}
  \label{kernelnorm}
  \|K_{w}\|_p= e^{\frac{1}{2}|w|^2}
  \end{align} for all  $w\in \CC$ and $0<p<\infty$.    This follows from a simple computation
   \begin{align*}
   \|K_{w}\|_p^p= \frac{p}{2\pi}\int_{\CC} e^{ p\Re(\overline{w}z)-\frac{p}{2}|z|^2}
 dA(z)
 =   e^{\frac{p}{2}|w|^2}\frac{p}{2\pi} \int_{\CC} e^{ \frac{p}{2}\big(2\Re(\overline{w}z)-|w|-|z|^2\big)}
 dA(z)\\
 = e^{\frac{p}{2}|w|^2}\bigg(\frac{p}{2\pi} \int_{\CC} e^{ -\frac{p}{2}|w-z|^2}
 dA(z)\bigg)= e^{\frac{p}{2}|w|^2}.
  \end{align*}
 It follows that  $k_w$ constitutes a unit norm sequence of functions  $\mathcal{F}_p$  and converges to zero in compact subset of $\CC$.
 This and the norms in \eqref{kernelnorm} will be used repeatedly in our further considerations.\\

In \cite{TM00,TM}, Mengestie studied the operators   $V_{(g,\psi)}$ and  $J_{(g,\psi)}$ on Fock spaces and characterized various properties of
the operators in terms of  the Berezin type integral transforms
\begin{align}
 B_{(|g|^p,\psi)} (w)&=\int_{\CC} \big|k_w(\psi(z))\big|^p
\frac{|g'(z)|^p e^{-\frac{p }{2}|z|^2}}{(1+|z|)^p} dA(z) \ \ \text{and}\nonumber\\
\widetilde B_{(|g|^p,\psi)} (w)&= \int_{\CC} \big|k_w(\psi(z))\big|^p
\frac{(|w|+1)^p|g(z)|^pe^{-\frac{p }{2}|z|^2}}{(1+|z|)^p}  dA(z).
\label{Berzin}
\end{align}
As will be indicated in Section~\ref{proof}, the  characterizations in \cite{TM00,TM} require to estimate
the $L^q, 0<q\leq \infty $  norms of the integral transforms $B_{(|g|^p,\psi)}$ and $\widetilde B_{(|g|^p,\psi)}$. Such types of characterizations have been also referred as  reproducing kernel thesis properties for the operators.  One of the main purposes of this work is to substantially improve these conditions and  provide
characterizations which are  rather simpler to apply.  Our new results will be  expressed in terms of  the functions
\begin{align*}
M_{(g,\psi)}(z)= \frac{|g'(z)|}{1+|z|}e^{\frac{1}{2}(|\psi(z)|^2-|z|^2)}\ \text{and} \ \widetilde M_{(g,\psi)}(z)= \frac{|g(z)|(1+|\psi(z)|)}{1+|z|}e^{\frac{1}{2}(|\psi(z)|^2-|z|^2)},
\end{align*} which are easier to handle  with than those class of integral  transforms in \eqref{Berzin}.
 \section{Bounded and Compact $V_{(g,\psi)}$ and  $J_{(g,\psi)}$ }
 In this section we characterized the boundedness and compactness properties of the operators $V_{(g,\psi)}$ and  $J_{(g,\psi)}$ in terms of the functions $M_{(g,\psi)}$ and $\widetilde M_{(g,\psi)}$. Our first main result reads as follows.
\begin{theorem}\label{bounded}
Let $0<p\leq q<\infty$ and $(g,\psi)$ be pairs of nonconstant  entire functions. Then
\begin{enumerate}
\item
$V_{(g,\psi)}: \mathcal{F}_p \to \mathcal{F}_q$ is  bounded if and only if
$ M_{(g,\psi)} \in L^\infty(\CC, dA)$.
\item  $V_{(g,\psi)}: \mathcal{F}_p \to \mathcal{F}_q$ is compact if and only if
 \begin{align*}
\lim_{|z|\to \infty}
M_{(g, \psi)}(z)= 0.\end{align*}
\item
$J_{(g,\psi)}: \mathcal{F}_p \to \mathcal{F}_q$ is  bounded if and only if
$ \widetilde M_{(g,\psi)}\in L^\infty(\CC, dA)$.
\item  $J_{(g,\psi)}: \mathcal{F}_p \to \mathcal{F}_q$ is   compact if and only if
 \begin{align*}\lim_{|z|\to \infty}
 \widetilde M_{(g, \psi)} (z)= 0.\end{align*}
\end{enumerate}
\end{theorem}
It should be remarked that  all the results above are independent of the Fock exponent $p$ in the domain spaces $\mathcal{F}_p$ apart from the condition that  it  is bounded by the exponent q in the target space $\mathcal{F}_q$. This is manifested  due to the fact that  the Fock spaces are nested in the sense that $\mathcal{F}_p \subseteq \mathcal{F}_q$ whenever $p\leq q$ (see Theorem 2.10, \cite{KZH1}). As one would expect, our results are different for the  cases $p\leq q$  and  $q < p$.  For the latter case, we have rather a stronger condition under which the boundedness implies  compactness as formulated  in our next main result.
\begin{theorem} \label{thsmall} Let $0<q<p<\infty$ and $(g,\psi)$ be pairs of nonconstant entire functions. Then
\begin{enumerate}
\item the following statements are equivalent.
\begin{enumerate}
 \item $V_{(g,\psi)}: \mathcal{F}_p \to \mathcal{F}_q$ is bounded;
 \item $V_{(g,\psi)}: \mathcal{F}_p \to \mathcal{F}_q$ is compact;
 \item  $M_{(g, \psi)}\in L^{\frac{pq}{p-q}}(\CC, dA)$.
 \end{enumerate}
 \item the following statements are also  equivalent.
 \begin{enumerate}
 \item  $J_{(g,\psi)}: \mathcal{F}_p \to \mathcal{F}_q$ is  bounded;
 \item $J_{(g,\psi)}: \mathcal{F}_p \to \mathcal{F}_q$ is  compact;
 \item $ \widetilde M_{(g, \psi)}\in L^{\frac{pq}{p-q}}(\CC, dA)$.
 \end{enumerate}
  \end{enumerate}
\end{theorem}
\section{ Schatten $\mathcal{S}_p$  class  membership of  $V_{(g,\psi)} $ }
The singular values  of a compact operator $T$ on a Hilbert space $\mathcal{H}$ are
the square roots of the positive eigenvalues of the operator $T^*T$, where
$T^*$ denotes the adjoint of $T$. Said differently, these are  simply the  positive
eigenvalues of $(T^*T)^{1/2}= |T|$.
 For $0<p<\infty$, the Schatten p-class of $\mathcal{H} $ which we denote it by $ \mathcal{S}_p(\mathcal{H})$ represents  the space of all compact operators $T$ on $H$ with  $\ell^p$ summable singular value sequences.  The space $\mathcal{S}_p(\mathcal{H})$ is a two sided ideal  in the algebra of all bounded linear operators on the space $\mathcal{H}$. We refer readers to the materials in  \cite{GO,KZH1} for a good overview of the subject.  \\
 Our next main result characterizes the Schatten $\mathcal{S}_p(\mathcal{F}_2)$ membership of the operators $V_{(g,\psi)}$ in terms of $L^p$ integrability condition for $ M_{(g,\psi)}$.
\begin{theorem}\label{Schatten}
Let $0<p<\infty$ and $(g,\psi)$ be pair of  entire functions and the map $V_{(g,\psi)}: \mathcal{F}_2\to \mathcal{F}_2$  be bounded. Then   $V_{(g,\psi)}$ belongs to the algebra $\mathcal{S}_p(\mathcal{F}_2)$ if and only if
 $ M_{(g,\psi)}\in L^{p}(\CC,dA)$.
 \end{theorem}
 The analogous problem  for $J_{(g,\psi)}$ has been already studied  in \cite{TM0} where  partly a  source of inspiration for this paper stems from.
 \section{Topological structures of  $V_{(g,\psi)} $ and $J_{(g,\psi)}$  } \label{topolotical}
 Over the past  few decades much effort has been paid  to study  various operator-theoretic
 properties of these operators  in terms of the function-theoretic
properties of the  inducing pairs $(g,\psi)$.
 On the contrary, there has been no effort devoted to understand the topological structures of the space of bounded  operators $V_{(g,\psi)} $ and $J_{(g,\psi)}$   equipped with  the operator norm topology. It is  worth noting that such structures have been studied for many other operators including  compact   differences  of  composition  operators.   We  refer  to  \cite{MO}  for
references and historical remarks in this topic.  In  this section, we use  our results from the proceeding sections and  begin  a study on such structures  and characterize the compact and Schatten  $\mathcal{S}_p$ class differences of the operators.
 \begin{theorem}\label{compactdifference} Let $0<p\leq q<\infty$, $(g_1,\psi_1)$ and $(g_2,\psi_2)$ be pairs of nonconstant entire functions,  and $V_{(g_1,\psi_1)}, \ V_{(g_2,\psi_2)}, \  J_{(g_1,\psi_1)}$ and $J_{(g_2,\psi_2)}$ be bounded operators from $\mathcal{F}_p$ into $\mathcal{F}_q .$ Then the operator
 \begin{enumerate}
 \item   $V_{(g_1,\psi_1)}-V_{(g_2,\psi_2)} $ is compact if and only if either both $V_{(g_1,\psi_1)}$ and $ V_{(g_2,\psi_2)} $ are compact or  $\psi_1= \psi_2 $ and
    \begin{align}
    \label{compdiffer}
    \lim_{|z|\to \infty}M_{(g_1-g_2,\psi_1)}(z)=0.\end{align}
 \item   $J_{(g_1,\psi_1)}-J_{(g_2,\psi_2)} $ is compact if and only if either both $J_{(g_1,\psi_1)}$ and $ J_{(g_2,\psi_2)} $ are compact or  $\psi_1= \psi_2 $ and
   \begin{align}
   \label{compdiffer1}
   \lim_{|z| \to \infty}\widetilde M_{(g_1-g_2,\psi_1)}(z)=0.\end{align}
 \end{enumerate}
  \end{theorem}
 An interesting future in the study of differences of operators  is the important factors arising from the cancellation property. From \eqref{compdiffer} and \eqref{compdiffer1}, such  factors  are $|\psi_1-\psi_2=0|$,  and $|g'_1-g_2'|$  and $|g_1-g_2|$.  The same impact is observed in the stronger  Schatten  $\mathcal{S}_p$ class  difference membership as stated in the next main result.
 \begin{theorem}\label{Schattendifference} Let $0<p<\infty$, $(g_1,\psi_1)$ and $(g_2,\psi_2)$ be pairs of nonconstant entire functions,  and $V_{(g_1,\psi_1)}, \ V_{(g_2,\psi_2)}, \  J_{(g_1,\psi_1)}$ and $J_{(g_2,\psi_2)}$ be bounded operators on $\mathcal{F}_2$. Then the operator
 \begin{enumerate}
 \item   $V_{(g_1,\psi_1)}-V_{(g_2,\psi_2)} $  belongs to the Schatten  $\mathcal{S}_p(\mathcal{F}_2)$ class if and only if either both $V_{(g_1,\psi_1)}$ and $ V_{(g_2,\psi_2)} $  belong to the    $\mathcal{S}_p(\mathcal{F}_2)$ class or  $\psi_1= \psi_2= \psi $ and
    \begin{align}
    \label{schattendiffer}
     \int_{\CC}  \frac{|g_1'(z)-g_2'(z)|^p}{(1+|z|)^p}e^{\frac{p}{2}(|\psi(z)|^2-|z|^2)} dA(z) <\infty.
  \end{align}
 \item   $J_{(g_1,\psi_1)}-J_{(g_2,\psi_2)} $ belongs to the Schatten  $\mathcal{S}_p(\mathcal{F}_2)$ class if and only if either both $J_{(g_1,\psi_1)}$ and $ J_{(g_2,\psi_2)} $  belong to the    $\mathcal{S}_p(\mathcal{F}_2)$ class or  $\psi_1= \psi_2 $ and
   \begin{align}
   \label{schattendiffer1}
  \int_{\CC}  \frac{|g_1(z)-g_2(z)|^p}{(1+|z|)^p}e^{\frac{p}{2}(|\psi(z)|^2-|z|^2)} dA(z) <\infty.\end{align}
 \end{enumerate}
  \end{theorem}
As noted  earlier,  if we  set $\psi(z)= z,$ then  the operators $ V_{(g,\psi)}$ becomes  the operators $V_g$.  By Theorem~\ref{compactdifference}, it turns out that a Volterra-type integral difference is compact if and only if both operators are compact. We record this as  part of the following  corollary since it is interests of its own.
 \begin{corollary}\label{cor1}
 \begin{enumerate}
 \item
 Let $0<p\leq q<\infty$,  $g_1\neq g_2$,  and $V_{g_1}$ and $ V_{g_2}$ be bounded operators from $\mathcal{F}_p$ into  $\mathcal{F}_q.$ Then  the  difference  operator
 \   $V_{g_1}-V_{g_2} $ is
 compact if and only if  both $V_{g_1}$ and $ V_{g_2} $ are compact.

 \item Let $g_1\neq g_2$ , $0<p<\infty$,  and  $V_{g_1}-V_{g_2} $  be compact operators on $\mathcal{F}_2$. Then the  operator
 $V_{g_1}-V_{g_2} $ belongs to the Schatten $S_p(\mathcal{F}_2)$  class if and only if $p>2.$
 \item Let $1\leq p<\infty$,   and $V_{g_1}$ and $V_{g_2} $  be bounded  operators on $\mathcal{F}_p$. That is $g_1(z)= a_1z^2+b_1z+c_1$ and $g_1(z)= a_2z^2+b_2z+c_2$. Then
the spectrum of  $V_{g_1}-V_{g_2}$ on  $\mathcal{F}_p$ is given by
\begin{align*}
\sigma \big(V_{g_1}-V_{g_2}\big)= \Big\{ \lambda \in \CC: |\lambda| \leq 2|a_1-a_2| \Big\}.
\end{align*}
 \end{enumerate}
 \end{corollary}
 This shows that  the difference of two   Volterra-type integral operators can not be nontrivially compact. A natural question to ask is whether  there exist compact difference of  multiplication operators $M_g$ and the other Volterra type operator $J_g$. Clearly, if any two of the operators in \eqref{parts} are  bounded, so is the  third one.
 The boundedness and compactness properties of  $M_g$ and  $J_g$ on Fock spaces were described in \cite{TM00}. The  result there  ensures that $M_g$ or $J_g$ is compact if and only
 if $g=0$. Thus, a difference of any two of such operators fails to be compact.
\section{Proof of the main results} \label{proof}
In this section we prove our results.
We may first give  a key  proposition that will be used repeatedly in our subsequent considerations.
\begin{proposition} \label{prop}
Let $(g,\psi)$ be  a pair of  nonconstant entire functions  and $0<p<\infty$. Then if any of
\begin{enumerate}
\item  $M_{(g,\psi)}, \ \widetilde M_{(g,\psi)},   B_{(|g|^p,\psi)}$ or $   \widetilde B_{(|g|^p,\psi)} $ belongs to $ L^\infty(\CC, dA)$,   then  $\psi(z)= az+b$ for some $|a|\leq 1$.
\item  $M_{(g,\psi)}(z), \ \widetilde M_{(g,\psi)}(z),   B_{(|g|^p,\psi)}(z), $ or $  \widetilde B_{(|g|^p,\psi)}(z) $  tends to zero as $|z| \to \infty,$ then $\psi(z)= az+b$ with $|a|<1$.
\end{enumerate}
 \end{proposition}
\emph{Proof.}
 The proof of the proposition for the parts  $M_{(g,\psi)}$ and $\widetilde M_{(g,\psi)}$   follows from a simple variant of the proof of Proposition 2.1 in \cite{TLE} or Lemma 2.3 of \cite{TM0}.
    On the other hand, for the integral transform  $B_{\psi}(|g|^p)$,  we have
    \vspace{-0.01in}
 \begin{align*}
  B_{\psi}(|g|^p)(w)= \int_{\CC} \frac{|k_w(\psi(z))|^p|g'(z)|^p}{(1+|z|)^pe^{\frac{p}{2}|z|^2}}  dA(z)
  \geq \int_{D(\zeta, 1)} \frac{|k_w(\psi(z))|^p|g'(z)|^p}{(1+|z|)^pe^{\frac{p}{2}|z|^2}}  dA(z)
  \end{align*} for all $w, \zeta \in \CC$, where  $D(\zeta, 1)$ is a disc of radius 1 and center $\zeta$.\\
    Since $1+z \footnote{The  notation $U(z)\lesssim V(z)$ (or
equivalently $V(z)\gtrsim U(z)$) means that there is a constant
$C$ such that $U(z)\leq CV(z)$ holds for all $z$ in the set of a
question. We write $U(z)\simeq V(z)$ if both $U(z)\lesssim V(z)$
and $V(z)\lesssim U(z)$.}\simeq 1+\zeta$ for all $z\in D(\zeta, 1)$, we further estimate the integral above  as
  \begin{align*}
\int_{D(\zeta, 1)} \frac{|k_w(\psi(z))|^p|g'(z)|^p}{(1+|z|)^pe^{\frac{p}{2}|z|^2}}  dA(z)\simeq \frac{1}{(1+|\zeta|)^p} \int_{D(\zeta, 1)} \frac{|k_w(\psi(z))|^p|g'(z)|^p}{(1+|z|)^pe^{\frac{p}{2}|z|^2}}  dA(z)\\
  \gtrsim \frac{1}{(1+|\zeta|)^p} |k_w(\psi(\zeta))g'(\zeta)|^p e^{-\frac{p}{2}|\zeta|^2}
 \end{align*} for all $\zeta\in \CC$.
Setting  $w= \psi(\zeta)$ in particular and applying \eqref{kernelnorm} gives
\begin{align}
 \label{series1}
  B_{\psi}(|g|^p)(\psi(\zeta))
  \gtrsim \frac{1}{(1+|\zeta|)^p} |k_{\psi(\zeta)}(\psi(\zeta))|^p |g'(\zeta)|^p e^{-\frac{p}{2}|\zeta|^2}\nonumber\\
  = \frac{|g'(\zeta)|^p}{(1+|\zeta|)^p}  e^{\frac{p}{2}|\psi(\zeta)|^2-\frac{p}{2}|\zeta|^2} = M_{(g,\psi)}^p(\zeta),
 \end{align} and the assertion follows  from the boundedness  condition on $M_{(g,\psi)}$.\\
 Similarly, if $\widetilde B_{(|g|^p,\psi)}$ is bounded, then
 \begin{align*}
  \widetilde B_{(|g|^p,\psi)} (w)= \int_{\CC} \big|k_w(\psi(z))\big|^p
\frac{(|w|+1)^p|g(z)|^pe^{-\frac{p }{2}|z|^2}}{(1+|z|)^p}  dA(z)\\
\geq \int_{D(\zeta, 1)} |k_w(\psi(z))|^p \frac{(1+|w|)^p|g(z)|^p}{(1+|z|)^p} e^{-\frac{p}{2}|z|^2} dA(z)\\
\gtrsim \frac{(1+|w|)^p}{(1+|\zeta|)^p} |k_w(\psi(\zeta))g(\zeta)|^p e^{-\frac{p}{2}|\zeta|^2},
\end{align*} and setting  $w= \psi(\zeta)$  and applying \eqref{kernelnorm} again results
\begin{align}
\label{series2}
\widetilde B_{(|g|^p,\psi)} (\psi(\zeta)) \gtrsim \frac{(1+|\psi(\zeta)|)^p}{(1+|\zeta|)^p} |k_{\psi(\zeta)}(\psi(\zeta))|^p |g(\zeta)|^p e^{-\frac{p}{2}|\zeta|^2}\nonumber\\
  = \frac{(1+|\psi(\zeta)|)^p|g(\zeta)|^p}{(1+|\zeta|)^p}  e^{\frac{p}{2}|\psi(\zeta)|^2-\frac{p}{2}|\zeta|^2} =\widetilde M_{(g,\psi)}^p(\zeta),
\end{align} and the assertion follows  from the boundedness  condition on $ \widetilde M_{(g,\psi)}$ and completes the proof.

 To prove our results, we also need the following crucial  Littlewood-Paley  type  estimate   for  Fock spaces \cite{Olivia}. The estimate asserts that for $0<p<\infty$ and $f\in \mathcal{F}_P$
 \begin{align}
 \label{Olivia}
 \int_{\CC} |f(z)|^p e^{-\frac{p}{2}|z|^2} dA(z) \simeq |f(0)|^p + \int_{\CC} \frac{|f'(z)|^p}{(1+|z|)^p} e^{-\frac{p}{2}|z|^2} dA(z).
 \end{align}
 We may now turn to the proofs of our main results.
  \subsection{Proof of Theorem~\ref{bounded}}
  By Theorem~1 of \cite{TM}, $V_{(g,\psi)}: \mathcal{F}_p \to \mathcal{F}_q$ is  bounded if and only if the  Berezin type integral transforms $B_{(|g|^q,\psi)}$ is bounded, and compact if and only if $B_{(|g|^q,\psi)}(z)  \to  0$ as $|z| \to \infty$. Seemingly, by Theorem~1 of \cite{TM00}, $J_{(g,\psi)}: \mathcal{F}_p \to \mathcal{F}_q$ is  bounded if and only if $\widetilde B_{(|g|^q,\psi)}$ is bounded and compact if and only if $\widetilde B_{(|g|^q,\psi)} \to 0$ as $|z| \to \infty$. Due to these, the proof of our theorem will be complete if  we prove the following two  results which we formulate them as lemmas as they  are interest of their own.
\begin{lemma}\label{lem1}
 Let  $(g,\psi)$ be a pair of nonconstant  entire functions on $\CC$. Then
 \begin{enumerate}
 \item  $ M_{(g,\psi)} \in L^\infty (\CC, dA)$ if and only if $ B_{(|g|^p,\psi)} \in L^\infty(\CC, dA)$ for some \\ $0<p<\infty$.
 \item $ \widetilde{M}_{(g,\psi)} \in L^\infty (\CC, dA)$ if and only if $  \widetilde B_{(|g|^p,\psi)}\in L^\infty(\CC, dA) $ for some\\  $0<p<\infty$.
  \end{enumerate}
  \end{lemma}
\emph{Proof.}
 (i). If the  Berezin type integral transform $B_{(|g|^p,\psi)}$ is bounded for some $0<p<\infty$, then the series of estimates  in \eqref{series1} implies that $ M_{(g,\psi)}$ is bounded. On the other hand, if $ M_{(g,\psi)}$ is bounded,  then by Proposition~\ref{prop},  we may set $\psi(z)= az+b$ and argue
 \begin{align*}
  B_{(|g|^p,\psi)}(w)= \int_{\CC}|k_w(\psi(z))|^p \frac{|g'(z)|^p}{(1+|z|)^p} e^{-\frac{p}{2}|z|^2} dA(z)\quad \quad \quad\quad \quad \quad \quad\quad \quad \quad\nonumber\\
  = \int_{\CC} |k_w(\psi(z))|^p  e^{-\frac{p}{2}|\psi(z)|^2} \bigg(\frac{|g'(z)|^p}{(1+|z|)^p} e^{\frac{p}{2}(|\psi(z)|^2-|z|^2)}\bigg) dA(z)\nonumber\\
  \leq \bigg(\sup_{z\in \CC} \frac{|g'(z)|^p}{(1+|z|)^p} e^{\frac{p}{2}(|\psi(z)|^2-|z|^2)} \bigg) \int_{\CC} |k_w(\psi(z))|^p  e^{-\frac{p}{2}|\psi(z)|^2} dA(z)\nonumber\\
  =  \bigg(  \sup_{z\in \CC} M_{(g,\psi)}^p(z)\bigg) \frac{1}{|a|^2}  \int_{\CC} |k_w(z)|^p  e^{-\frac{p}{2}|z|^2} dA(z)  = \frac{1}{|a|^2} \sup_{z\in\CC} M_{(g,\psi)}^p (z),
 \end{align*} where the last equality follows by \eqref{kernelnorm} and hence the claim. \\
 Similarly, one side of the statement in part (ii) follows from the series of the estimations made leading to \eqref{series2}. On the other hand, applying Proposition~\ref{prop} and \eqref{kernelnorm} again,
 \begin{align*}
  \widetilde{B}_{(|g|^p,\psi)}(w)  = \int_{\CC} \frac{|k_w(\psi(z))|^p }{ e^{\frac{p}{2}|\psi(z)|^2}} \bigg(\frac{(1+|w|)^p(1+|\psi(z)|)^p|g(z)|^p}{(1+|\psi(z)|)^p(1+|z|)^p} e^{\frac{p}{2}(|\psi(z)|^2-|z|^2)}\bigg) dA(z)\nonumber\\
  \lesssim \bigg(\sup_{z\in \CC} \frac{(1+|\psi(z)|)^p|g(z)|^p}{(1+|z|)^p} e^{\frac{p}{2}(|\psi(z)|^2-|z|^2)} \bigg) \int_{\CC} |k_w'(\psi(z))|^p  e^{-\frac{p}{2}|\psi(z)|^2} dA(z)\nonumber\\
  =  \bigg(  \sup_{z\in \CC} \widetilde{M}_{(g,\psi)}^p(z)\bigg) \frac{1}{|a|^2}  \int_{\CC}\frac{ |k_w'(z)|^p}{(1+|\psi(z)|)^p}  e^{-\frac{p}{2}|z|^2} dA(z).
 \end{align*}
 Invoking \eqref{Olivia}, \eqref{kernelnorm} and Proposition~\ref{prop} we estimate the above last integral by
 \begin{align*}
  \int_{\CC}\frac{ |k_w'(\psi(z))|^p}{(1+|\psi(z)|)^p}  e^{-\frac{p}{2}|\psi(z)|^2} dA(z)\simeq \|k_w\|_p^p = 1,
 \end{align*} from which we have that
 \begin{align*}
 \widetilde{B}_{(|g|^p,\psi)}(w) \lesssim \sup_{z\in\CC} \widetilde{M}_{(g,\psi)}^p(z)  \|k_w\|_p^p=  \sup_{z\in\CC} \widetilde{M}_{(g,\psi)}^p(z),
 \end{align*} and completes the proof.
  \begin{lemma}\label{lemcomp}
  Let $(g,\psi)$ be a pair of  nonconstant entire functions on $\CC$. Then
  \begin{enumerate}
  \item
  $ M_{(g,\psi)}(w)\to 0 $ as $|w| \to \infty$  if and only if $  B_{(|g|^p,\psi)}(w)  \to  0 $ as $|w| \to \infty$ for some $0<p<\infty$.
  \item $ \widetilde {M}_{(g,\psi)}(w)\to 0$  as $|w| \to \infty$  if and only if $   \widetilde B_{(|g|^p,\psi)}(w)  \to  0 $ as $|w| \to \infty$ for some $0<p<\infty$.
  \end{enumerate}
 \end{lemma}
\emph{Proof.}
(i) One side of the statement follows easily from the estimates in \eqref{series1}. We shall proceed to show the other side, and assume that  $ M_{(g,\psi)} (w)\to 0 $ as $|w| \to \infty$.

  Then we estimate
 \begin{align*}
   B_{(|g|^p,\psi)}(w)= \int_{\CC} \frac{ |k_w(\psi(z))|^p|g'(z)|^p}{(1+|z|)^p e^{\frac{p}{2}|z|^2}}  dA(z)  = \int_{|z|\leq| w|} \frac{|k_w(\psi(z))|^p}{e^{\frac{p}{2}|\psi(z)|^2}}   M_{(g,\psi)}^p(z) dA(z)\nonumber\\
  + \int_{|z|> |w|} \frac{|k_w(\psi(z))|^p M_{(g,\psi)}^p(z)}{  e^{\frac{p}{2}|\psi(z)|^2}} dA(z)
    \lesssim  \bigg(\sup_{z: |z|\leq |w|} |k_w(\psi(z))|^p\bigg) \int_{|z|\leq |w|} \frac{ M_{(g,\psi)}^p(z)}{  e^{\frac{p}{2}|\psi(z)|^2}} dA(z)\\
    + \bigg(\sup_{z: |z|> |w|} M_{(g,\psi)}^p(z)\bigg)  \int_{|z|> |w|} |k_w(\psi(z))|^p  e^{-\frac{p}{2}|\psi(z)|^2} dA(z)\\
    \end{align*}
    By  Proposition~\ref{prop} and  the assumption that $\psi$ is nonconstant, it follows
    \begin{align*}
    \int_{|z|\leq |w|} M_{(g,\psi)}^p(z) e^{-\frac{p}{2}|\psi(z)|^2}dA(z) \leq \bigg(\sup_{z\in \CC} M_{(g,\psi)}^p(z)\bigg) \int_{|z|\leq |w|}e^{-\frac{p}{2}|\psi(z)|^2}dA(z) <\infty
    \end{align*}
    from which we deduce
    \begin{align*}
   B_{(|g|^p,\psi)}(w)
    \lesssim \sup_{z: |z|\leq |w|} |k_w(\psi(z))|^p +
   \sup_{z: |z|> |w|} M_{(g,\psi)}^p(z) \quad\quad \quad \quad \quad \quad \quad\quad \quad \quad \\
    = \sup_{z: |z|\leq |w|} |k_w(\psi(z))|^p +
    \sup_{z: |z|> |w|} M_{(g,\psi)}^p(z).
     \end{align*}
     The first term in the last sum converges to zero as $|w|\to \infty$ since the normalized reproducing kernel converges to zero on compact subsets. On the other hand, by hypothesis, we have
       \begin{align*}
       \sup_{z: |z|> |w|} M_{(g,\psi)}(z) \to 0 \ \ \text{as}\  |w| \to \infty.
       \end{align*} and hence  the assertion.

     (ii) From the inequalities in \eqref{series2}, we easily deduce one of the implications for this part as well. On the other hand, if $ \widetilde {M}_{(g,\psi)}(w)\to 0$  as $|w| \to \infty$, then  arguing as above and eventually applying \eqref{Olivia} and \eqref{kernelnorm}
     \begin{align*}
 \widetilde B_{(|g|^p,\psi)}(w)= \int_{\CC}|k_w(\psi(z))|^p \frac{(1+|w|)^p|g(z)|^p}{(1+|z|)^p} e^{-\frac{p}{2}|z|^2} dA(z)\quad \quad \quad \quad \quad \quad \quad  \nonumber\\
   \simeq  \int_{|z|\leq| w|} \frac{|k'_w(\psi(z))|^p}{(1+|\psi(z)|)^p}  e^{-\frac{p}{2}|\psi(z)|^2} \widetilde {M}_{(g,\psi)}^p(z) dA(z)\\
   + \int_{|z|> |w|} \frac{|k_w'(\psi(z))|^p}{(1+|\psi(z)|)^p}  e^{-\frac{p}{2}|\psi(z)|^2} \widetilde {M}_{(g,\psi)}^p(z) dA(z)\nonumber\\
    \lesssim  \bigg(\sup_{z: |z|\leq |w|} |k'_w(\psi(z))|^p\bigg) \int_{|z|\leq |w|} \widetilde {M}_{(g,\psi)}^p(z)e^{-\frac{p}{2}|\psi(z)|^2}dA(z)\\
    + \bigg(\sup_{z: |z|> |w|} \widetilde {M}_{(g,\psi)}^p(z)\bigg)  \int_{|z|> |w|}\frac{ |k_w'(\psi(z))|^p}{(1+|\psi(z)|)^p}  e^{-\frac{p}{2}|\psi(z)|^2} dA(z)\\
    \lesssim \sup_{z: |z|\leq |w|} |k_w'(\psi(z))|^p +
    \bigg(\sup_{z: |z|> |w|} \widetilde {M}_{(g,\psi)}^p(z)\bigg)  \|k_w\|_p^p
     \end{align*} from which the assertion follows as before since  $k_w'$ converges uniformly
     on compact subsets and  $\sup_{z: |z|> |w|}\widetilde {M}_{(g,\psi)}(z)\to  0$  when $|w|  \to \infty$.
       \subsection{Proof of Theorem~\ref{thsmall}}
  By Theorem~2 of \cite{TM}, $V_{(g,\psi)}: \mathcal{F}_p \to \mathcal{F}_q$ is  bounded (compact) if and only if the  integral transform $B_{(|g|^q,\psi)}$ belongs to $L^{\frac{p}{p-q}}(\CC,dA)$. By Theorem~2 of \cite{TM00}, the operator  $J_{(g,\psi)}: \mathcal{F}_p \to \mathcal{F}_q$ is also  bounded (compact) if and only if $\widetilde B_{(|g|^q,\psi)}$ belongs to $L^{\frac{p}{p-q}}(\CC,dA)$. Thus, our theorem will be proved  once we prove the following key  lemma which is  again interest of its  own.
 \begin{lemma}\label{lem2}
  Let $0<q<p<\infty$ and $(g,\psi)$ be pairs of nonconstant  entire functions on $\CC$. Then
  \begin{enumerate}
 \item $ M_{(g,\psi)} \in L^\frac{pq}{p-q}(\CC, dA)$ if and only if $  B_{(|g|^q,\psi)} \in L^\frac{p}{p-q}(\CC, dA).$
 \item $ \widetilde {M}_{(g,\psi)} \in L^\frac{pq}{p-q}(\CC, dA)$ if and only if $  \widetilde {B}_{(|g|^q,\psi)} \in L^\frac{p}{p-q}(\CC, dA).$
 \end{enumerate}
 \end{lemma}
\emph{Proof.}
  (i). Making use of the estimate in \eqref{series1} and Proposition~\ref{prop} namely that  $\psi(z)= az+b$, we estimate
 \begin{align*}
\int_{\CC} M_{(g,\psi)}^{\frac{pq}{p-q}}(z) dA(z) \lesssim \int_{\CC}   B_{(|g|^q,\psi)}^{\frac{p}{p-q}}(\psi(z)) dA(z)=\frac{1}{|a|^2} \int_{\CC}  B_{(|g|^q,\psi)}^{\frac{p}{p-q}}(z) dA(z),
  \end{align*} and the necessity of the conditions follows.

  To prove the sufficiency, we may consider
  \begin{align*}
  \int_{\CC}    B_{(|g|^q,\psi)}^{\frac{p}{p-q}}(w)dA(w)= \int_{\CC}\bigg( \int_{\CC} |k_w(\psi(z))|^q \frac{|g'(z)|^q}{(1+|z|)^q}e^{-\frac{q}{2}|z|^2} dm(z)\bigg)^{\frac{p}{p-q}} dA(w).
    \end{align*}
  Applying H\"older's inequality to the inner integral above gives
    \begin{align*}
    H(w):= \bigg( \int_{\CC} |k_w(\psi(z))|^q \frac{|g'(z)|^q}{(1+|z|)^q}e^{-\frac{q}{2}|z|^2}dA(z)\bigg)^{\frac{p}{p-q}}\quad \quad \quad \quad \quad \quad \quad \quad\\
    \leq \int_{\CC} |k_w(\psi(z))|^q \frac{|g'(z)|^{\frac{qp}{p-q}}}{(1+|z|){\frac{qp}{p-q}}}e^{-\frac{qp}{2(p-q)}|z|^2} e^{\frac{q}{2} \big(\frac{p}{p-q}-1\big)|\psi(z)|^2} dA(z)\nonumber\\
    \times \bigg( \int_{\CC} |k_w(\psi(z))|^q e^{-\frac{q}{2}|\psi(z)|^2} dA(z)\bigg)^{\frac{q}{p-q}} \quad \quad \nonumber\\
    \lesssim \int_{\CC} |k_w(\psi(z))|^q \bigg|\frac{g'(z)}{1+|z|}\bigg|^{\frac{qp}{p-q}} e^{-\frac{qp}{2(p-q)}|z|^2} e^{\frac{q}{2} \big(\frac{p}{p-q}-1\big)|\psi(z)|^2} dA(z).
    \end{align*}
    Making use of  Fubini's Theorem and \eqref{kernelnorm}, we  further  estimate
    \begin{align*}
  \int_{\CC}   B_{(|g|^q,\psi)}^{\frac{p}{p-q}}(w)dA(w)= \int_{\CC}H(w) dA(w) \quad \quad \quad \quad \quad \quad \quad \quad \nonumber\\
  = \int_{\CC}\bigg|\frac{g'(z)}{1+|z|}\bigg|^{\frac{qp}{p-q}}e^{-\frac{qp}{2(p-q)}|z|^2} e^{\frac{q}{2} \big(\frac{p}{p-q}-1\big)|\psi(z)|^2} \int_{\CC} |K_{\psi(z)}(w)|^q e^{-\frac{q}{2}|w|^2} dA(w)dA(z)\nonumber\\
  \simeq \int_{\CC} \bigg|\frac{g'(z)}{(1+|z|)}\bigg|^{\frac{qp}{p-q}}e^{\frac{qp}{2(p-q)}\big( |\psi(z)|^2-|z|^2\big)} dA(z)
  =  \int_{\CC} M_{(g,\psi)}^{\frac{pq}{p-q}}(z) dA(z)<\infty.
    \end{align*}
    (ii). The necessity of the condition follows from the inequalities in \eqref{series2} and the fact that
$\psi(z)= az+b$ as ensured by Proposition~\ref{prop}.

 To prove the sufficiency,  we may set
\begin{align*}
T(w):= \bigg( \int_{\CC} |k_w(\psi(z))|^q \frac{(1+|w|)^q|g(z)|^q}{(1+|z|)^q}e^{-\frac{q}{2}|z|^2}dA(z)\bigg)^{\frac{p}{p-q}},
\end{align*}
 and   write
\begin{align*}
 \int_{\CC} \widetilde{ B}_{(|g|^q,\psi)}^{\frac{p}{p-q}}(w)dA(w)= \int_{\CC}T(w) dA(w).
\end{align*} Apply  H\"older's inequality  and subsequently  Fubini's Theorem, \eqref{Olivia},  and  \eqref{kernelnorm}  we  observe that   the above integral is bounded by
\begin{align*}
\int_{\CC}\bigg|\frac{(1+|\psi(z)|)g(z)}{1+|z|}\bigg|^{\frac{qp}{p-q}}e^{-\frac{qp}{2(p-q)}|z|^2} e^{\frac{q}{2} \big(\frac{p}{p-q}-1\big)|\psi(z)|^2}\quad \quad \quad \quad \quad \quad \quad \quad \quad \quad\nonumber\\
\times \int_{\CC} \frac{|k_{\psi(z)}'(w)|^q}{(1+|\psi(w)|)^q} e^{-\frac{q}{2}|w|^2} dA(w)dA(z)\quad \quad \nonumber\\
  \simeq \int_{\CC} \bigg|\frac{(1+|\psi(z)|)g(z)}{1+|z|}\bigg|^{\frac{qp}{p-q}}e^{\frac{qp}{2(p-q)}\big( |\psi(z)|^2-|z|^2\big)} dA(z)
  =  \int_{\CC} \widetilde{M}_{(g,\psi)}^{\frac{pq}{p-q}}(z) dA(z),
\end{align*} and completes the proof.
   \subsection{Proof of Theorem~\ref{Schatten}}
  By Theorem~4 of \cite{TM},  the operator $V_{(g,\psi)}: \mathcal{F}_2 \to \mathcal{F}_2$ belongs to the Schatten $S_p(\mathcal{F}_2)$ class  if and only  if the Berezin type integral transform $B_{(|g|^2,\psi)}$ belongs to $L^{\frac{p}{2}}(\CC,dA)$. Thus,  we only need to establish the following lemma.
 \begin{lemma}
 \label{lem3} Let $0<p<\infty$ and $(g,\psi)$ be a  pair of  nonconstant  entire functions on $\CC$. Then
  $ M_{(g,\psi)} \in L^p(\CC, dA)$ if and only if $ B_{(|g|^2,\psi)} \in L^\frac{p}{2}(\CC, dA).$
 \end{lemma}
 We remark in passing that  Lemma~\ref{lem3} does not follow from Lemma~\ref{lem2}.

 \emph{Proof.}
 Applying the inequalities in  \eqref{series1} and Proposition~\ref{prop} we have
 \begin{align*}
\int_{\CC} M_{(g,\psi)}^p(z) dA(z) \leq \int_{\CC} B_{(|g|^2,\psi)}^{\frac{p}{2}}(\psi(z)) dA(z)
=\frac{1}{|a|^2} \int_{\CC}  B_{(|g|^2,\psi)}^{\frac{p}{2}}(z) dA(z),   \end{align*}   from which  one side of the  assertion in the lemma  follows.\\
 For the remaining part, we consider two different cases.

 \emph{Case 1}. If  $0<p<2,$ then using  \eqref{Olivia} and  the fact that $\mathcal{F}_p \subset \mathcal{F}_2$ for $0<p\leq 2$,
 \begin{align}
\int_{\CC} |k_w(\psi(\zeta))|^2 \frac{|g'(\zeta)|^2e^{-|\zeta|^2}}{(1+|\zeta|)^2}  dA(\zeta) \qquad \qquad  \qquad \qquad \qquad  \qquad \qquad \qquad  \nonumber\\
\simeq \int_{\CC}\Bigg|\int_{0}^z k_w(\psi(\zeta))g'(\zeta) dA(\zeta)\Bigg|^2 e^{-|z|^2} dA(z)\nonumber\\
\lesssim \Bigg( \int_{\CC}\bigg|\int_{0}^z k_w(\psi(\zeta))g'(\zeta) dA(\zeta)\bigg|^p e^{-\frac{ p}{2}|z|^2} dA(z)\Bigg)^{\frac{2}{p}}\nonumber\\
\simeq \Bigg(\int_{\CC} |k_w(\psi(\zeta))|^p \frac{|g'(\zeta)|^pe^{-\frac{ p}{2} |\zeta|^2}}{(1+|\zeta|)^p}  dA(\zeta)\Bigg)^{\frac{2}{p}}.\nonumber
\end{align}
Using  this, Fubini's Theorem and \eqref{kernelnorm},  we further estimate
\begin{align*}
\int_{\CC} B_{(|g|^2,\psi)}^{\frac{p}{2}}(w) dA(w)\lesssim  \int_{\CC}\int_{\CC} |k_w(\psi(\zeta))|^p \frac{|g(\zeta)|^pe^{-\frac{ p}{2} |\zeta|^2}}{(1+|\zeta|)^p}  dA(\zeta) dA(w)\\
\int_{\CC} \frac{|g'(\zeta)|^pe^{-\frac{ p}{2} |\zeta|^2}}{(1+|\zeta|)^p} \int_{\CC} |K_{\psi(\zeta)}(w)|^pe^{-\frac{p}{2}|w|^2} dA(w)dA(\zeta)\\
\simeq \int_{\CC} \frac{|g'(\zeta)|^pe^{\frac{ p}{2}(\psi(\zeta)- |\zeta|^2}}{(1+|\zeta|)^p} dA(\zeta)
= \int_{\CC} M_{(g,\psi)}^p(\zeta) dA(\zeta)<\infty.
\end{align*}
\emph{Case 2}. If $p>2$, then applying  H\"older's inequality  we get
\begin{align*}
S(w): = \Bigg( \int_{\CC} |k_w(\psi(z))|^2 \frac{|g'(z)|^2e^{-|z|^2}}{ (1+|z|)^2} dA(z)\Bigg)^{\frac{p}{2}} \quad \quad \quad \quad \quad  \quad \quad \quad \quad \quad  \quad \quad \nonumber\\
\leq \Bigg( \int_{\CC} | k_w(\psi(z))|^2 \frac{|g'(z)|^p e^{- \frac{p}{2}|z|^2}}{ (1+|z|)^p} e^{ (\frac{p}{2}-1) |\psi(z)|^2}dA(z)\Bigg)\nonumber\\
 \times \ \  \Bigg(\int_{\CC}| k_w(\psi(z))|^2 e^{-|\psi(z)|^2}  dA(z)\Bigg) ^{\frac{p-2}{2}}.
    \end{align*}
  Making a change of variables and applying \eqref{kernelnorm}  again yields
    \begin{align*}
\int_{\CC}| k_w(\psi(z))|^2 e^{-|\psi(z)|^2}  dA(z)\simeq 1,
    \end{align*}
     which implies that
    \begin{align*}
S(w) \lesssim \int_{\CC} | k_w(\psi(z))|^2 \frac{|g'(z)|^p e^{- \frac{p}{2}|z|^2}}{ (1+|z|)^p} e^{(\frac{p}{2}-1) |\psi(z)|^2}dA(z).
    \end{align*} From this, Fubini's Theorem  and \eqref{kernelnorm}   we have the estimate
       \begin{align*}
\int_{\CC} B_{(|g|^2,\psi)}^{\frac{p}{2}}(w) dA(w) = \int_{\CC} S(w) dA(w)\quad \quad \quad \quad \quad \quad\quad \quad \quad\quad \quad \quad\quad \\
 \lesssim \int_{\CC} \int_{\CC} | k_w(\psi(z))|^2 \frac{|g'(z)|^p e^{- \frac{p}{2}|z|^2}}{ (1+|z|)^p} e^{(\frac{p}{2}-1) |\psi(z)|^2}dA(z) dA(w).\nonumber\\
= \int_{\CC} \frac{|g'(z)|^p e^{- \frac{p}{2}|z|^2}}{ (1+|z|)^p} e^{(\frac{p}{2}-1) |\psi(z)|^2} \int_{\CC} |K_{\psi(z)}(w)|^2 e^{-|w|^2}dA(w) dA(z)\nonumber\\
\simeq  \int_{\CC} \frac{|g'(z)|^p }{ (1+|z|)^p} e^{\frac{p}{2}( |\psi(z)|^2-|z|^2)} dA(z)
= \int_{\CC} M_{(g,\psi)}^p(z) dA(z)<\infty.
  \end{align*}
 \subsection{Proof of Theorem~\ref{compactdifference}} The sufficiency of the condition in the  theorem follows  from Theorem~\ref{bounded}. Thus, we may assume that the difference  $V_{(g_1,\psi_1)}-V_{(g_2,\psi_2)} $ is compact and proceed to show the necessity. If one of either $V_{(g_1,\psi_1)} $ or $V_{(g_2,\psi_2)} $ is compact, so is the other one which follows from the algebra of compact operators.  It follows that either both operators are compact or both are noncompact.  Thus, we may assume the later. To this end,  by Theorem~\ref{bounded}, there exist two positive numbers $\alpha_1$ and $\alpha_2$ such that
 \begin{align*}
 \alpha_1= \limsup_{|z|\to \infty} M_{(g_1,\psi_1)}(z) \ \ \text{and}
 \end{align*}
 \begin{align*}
 \alpha_2= \limsup_{|z|\to \infty} M_{(g_2,\psi_2)}(z).
 \end{align*}
 Since $k_w$ is weakly convergent, compactness of  $V_{(g_1,\psi_1)}-V_{(g_2,\psi_2)} $ implies that
 \begin{align}
 \label{k}
 \| \big(V_{(g_1,\psi_1)}-V_{(g_2,\psi_2)}\big)k_w\|_q \to 0 \quad \text{as} \quad \ |w|\to \infty.
 \end{align}
 On the other hand, using \eqref{Olivia} and the techniques leading to \eqref{series1}, we have
 \begin{align}
 \label{kk}
\| \big(V_{(g_1,\psi_1)}-V_{(g_2,\psi_2)}\big)k_w\|_q \simeq\bigg( \int_{\CC} \frac{|g_1'(z)k_w(\psi_1(z))-g_2'(z)k_w(\psi_2(z))|^q}{(1+|z|)^q e^{\frac{q}{2}|z|^2}}
dA(z)\bigg)^{\frac{1}{q}}\nonumber\\
\gtrsim  \frac{1}{1+|z|} \big|g_1'(z)k_{\psi_1(z)}(\psi_1(z))-g_2'(z)k_{\psi_1(z)}(\psi_2(\psi_1))\big| e^{-\frac{1}{2}|z|^2}
 \end{align} after setting $w= \psi_1(z)$.The right-hand side above can be estimated further as
 \begin{align}
 \label{kkk}
\frac{1}{1+|z|} \big|g_1'(z)k_{\psi_1(z)}(\psi_1(z))-g_2'(z)k_{\psi_1(z)}(\psi_2(\psi_1))\big| e^{-\frac{1}{2}|z|^2}\quad \quad \quad \quad \quad \quad \nonumber\\
 \geq \big|M_{(g_1,\psi_1)}(z)-M_{(g_2,\psi_2)}(z)e^{-\frac{1}{2}|\psi_1(z)-\psi_2(z)|^2}\big|.
 \end{align} Since  $V_{(g_1,\psi_1)}$ and $V_{(g_2,\psi_2)}$ are bounded operators,  by Theorem~\ref{bounded} and Proposition~\ref{prop},  $\psi_1$ and $\psi_2$ have  the linear forms  $\psi_1(z)= az+b$ and $\psi_2(z)= cz+d$ where  $|a|\leq 1$ and $|c|\leq 1.$ It follows from this that if $a\neq c$, then
 \begin{align*}
 \lim_{|z| \to \infty}e^{-\frac{1}{2}|\psi_1(z)-\psi_2(z)|^2}=0.
 \end{align*} Combining this with \eqref{k}, \eqref{kk}, \eqref{kkk} and the triangle inequality, we deduce
 \begin{align*}
 M_{(g_1,\psi_1)}(z) \leq \frac{M_{(g_2,\psi_2)}(z)}{e^{\frac{1}{2}|\psi_1(z)+\psi_2(z)|^2}}+ \Big|M_{(g_1,\psi_1)}(z)- \frac{M_{(g_2,\psi_2)}(z)}{e^{\frac{1}{2}|\psi_1(z)+\psi_2(z)|^2}}\Big|
 \end{align*} which implies
 \begin{align*}
 \lim_{|z| \to \infty}M_{(g_1,\psi_1)}(z) =0.
 \end{align*} It follows from this and Theorem~\ref{bounded} that $V_{(g_1,\psi_1)}$ is compact which contradicts our assumption.  Thus, we must have $a= c$. Taking this into account,
 \begin{align*}
 \alpha_1-\alpha_2 e^{-\frac{1}{2}|b-d|^2} \leq \limsup_{|z|\to \infty} \big|M_{(g_1,\psi_1)}(z)-M_{(g_2,\psi_2)}(z)e^{-\frac{1}{2}|b-d|^2}\big|\nonumber\\
 \lesssim \limsup_{|z|\to \infty}  \| \big(V_{(g_1,\psi_1)}-V_{(g_2,\psi_2)}\big)k_{\psi_1(z)}\|_q=0
 \end{align*} from which we get
 \begin{align}
 \label{kkkk}
  \alpha_1\leq \alpha_2 e^{-\frac{1}{2}|b-d|^2}.
 \end{align}On the other hand, if we repeat the above process by setting $w= \psi_2(z)$, we  get
 \begin{align*}
  \alpha_2\leq \alpha_1 e^{-\frac{1}{2}|b-d|^2}.
 \end{align*}
 From  this and \eqref{kkkk}, we find
 \begin{align*}
\alpha_1\leq \alpha_2 e^{-\frac{1}{2}|b-d|^2} \leq \alpha_1 e^{-|b-d|^2}\leq \alpha_1,
 \end{align*} which holds only if $b=d$. This shows  that $\psi_1= \psi_2 $ and hence the necessity of the condition follows from Theorem~\ref{bounded}.\\
  The proof of part (ii) follows in   a similar fashion.
   \subsection{Proof of Theorem~\ref{Schattendifference}}
   Since all Schatten $S_p(\mathcal{F}_2)$ class operators are compact, we can assume that $V_{(g_1, \psi_1)}- V_{(g_2, \psi_2)}$ is compact. Then by Theorem~\ref{compactdifference},  the difference is compact if and only if either both $V_{(g_1, \psi_1)}$ and $ V_{(g_2, \psi_2)}$ are compact or $ \psi_1= \psi_2= \psi$  and condition \eqref{compdiffer} holds. If  both are compact and the difference is in  the $\mathcal{S}_p$ class, then either both are in the $S_p(\mathcal{F}_2)$ class or both are not. In the latter case, following the same argument as in the proof of Theorem~\ref{compactdifference}, the  assumption that $V_{(g_1, \psi_1)}- V_{(g_2, \psi_2)}$ belongs to $\mathcal{S}_p$ implies $\psi_1= \psi_2= \psi.$
   On the other hand, since $V_{(g_1, \psi_1)}- V_{(g_2, \psi_2)}= V_{(g_1-g_2, \psi)}$ is itself a generalized Volterra-type integral operator induced by the pair of symbols $(g_1-g_2, \psi)$, by Theorem~\ref{Schatten},  it belongs to the $\mathcal{S}_p$ class if and only $ M_{(g_1-g_2,\psi)}\in L^{p}(\CC,dA)$.  That is
   \begin{align*}
   \int_{\CC} \frac{|g'(z)|^p}{(1+|z|)^p}e^{\frac{p}{2}(|\psi(z)|^2-|z|^2)} dA(z) <\infty
   \end{align*} which completes the proof of part (i) in the theorem.

   The proof of part (ii) follows in a similar manner.
       \subsection{Proof of Corollary~\ref{cor1}}
  Part (i) of the corollary follows once  by setting $\psi_1(z)= \psi_2(z)= z$ in Theorem~\ref{compactdifference}. Thus, we shall verify part (ii).  By part (i),  if the difference $V_{g_1}-V_{g_2}$ is compact,  then both $V_{g_1}$ and $V_{g_2}$ are compact and hence  Corollary~2 of \cite{TM} gives the representation    $g_1(z)= az+b$ and $g_2(z)= cz+d$.  Since the integral is linear, we also have the relation
  \begin{align*}
  V_{g_1}-V_{g_2}= V_{g_1-g_2}.
  \end{align*}
  On the other hand,  $V_{g_1-g_2} = V_{(g_1-g_2, z)}$. Then by Theorem~\ref{Schatten}, the difference operator  $V_{g_1-g_2}$ belongs to the Schatten $S_p$ class if and only if
  \begin{align}
  \label{simple}
  \int_{\CC} \frac{|(g_1'-g_2')(z)|^p}{(1+|z|)^p} dA(z)= \int_{\CC} \frac{|a+b-c-d|^p}{(1+|z|)^p}dA(z)<\infty.
  \end{align}
 Since $g_1\neq g_2$, we easily see that \eqref{simple} holds  only if $p>2$.

 Part (iii). By linearity of the integral, $V_{g_1}-V_{g_2} = V_{g_1-g_2}=  V_{g_3}$ where \begin{align*}
 g_3(z)= (a_1-a_2) z^2 + (b_1-b_2) z+ c_1+c_2.\end{align*} Then by  Theorem~1.3 of \cite{TM000},
 \begin{align*}
 \sigma(V_{g_3})= \{0\} \cup \overline{\Big \{  \lambda \in \CC \setminus \{ 0 \}:  e^{\frac{g_3}{\lambda}}\in \mathcal{F}_p \Big\}} \quad \quad \quad \quad \quad \quad  \quad \quad \quad  \quad \quad \quad  \\
 = \Big\{ \lambda \in \CC: |\lambda| \leq 2 |a_1-a_2|\Big\}=  \sigma\big(V_{g_1}-V_{g_2}\big).
 \end{align*}

\end{document}